\documentclass{article}
\usepackage{amsfonts}
\usepackage{amsmath}
\usepackage{amssymb}

\newtheorem{theorem}{Theorem}[section]

\newtheorem{lemma}[theorem]{Lemma}

\numberwithin{equation}{section}

\newcommand\func{\operatorname}

\begin{document}

\title{Some curvature restrictions on $3$-dimensional generalized ($\kappa
,\mu $)-contact metric manifolds}
\author{$^{1}$Manoj Ray Bakshi, $^{2}$Kanak Kanti Baishya}
\date{}
\maketitle

\begin{abstract}
The object of the present study is to study 3-dimensional
generalized ($\kappa ,\mu $)-contact metric manifolds with $\tilde{W}\cdot
R=0$ and $\tilde{W}\cdot H=0$ to cover all the eight equivalent classes
given in \cite{Shaikh2}.
\end{abstract}

\footnotetext{%
}

\section{Introduction}

In 1995, Blair, Koufogiorgos and Papantoniou \cite{Blairkofp}, introduced
the concept of ($\kappa ,\mu $)-contact metric manifold with several
examples. Thereafter, E. Boeckx \cite{Boec} classified such a manifold.
Considering $\kappa ,\mu $ as smooth functions, Koufogiorgos and
Tschlias \cite{Koufogiorgos1} defined the notion of generalized
($\kappa ,\mu $)-contact metric manifolds and claimed its existence for
3-dimensional case whereas, such a manifold does not exist for dimension
greater than 3. The 3-dimensional generalized ($\kappa
,\mu $)-contact metric manifolds are also studied in \cite{Goou},
\cite{Karatsobanis}, \cite{Koufogiorgos1}, \cite{Koufogiorgos2}, \cite%
{Koufogiorgos3} and \cite{Shaikh}.

Recently, in \cite{Baishya}, the authors have introduced and studied
generalized quasi-conformal curvature tensor (which will be denoted by $%
\tilde{W}$-tensor) in the frame of $N(k,\mu )$-manifold. For a ($2n+1)$%
-dimensional manifold, $\tilde{W}$-tensor is defined as follows%
\begin{eqnarray}
\tilde{W}(X,Y)Z&=&R(X,Y)Z +\alpha [S(Y,Z)X-S(X,Z)Y]\notag \\
&+&\beta [g(Y,Z)QX-g(X,Z)QY]
\notag \\
&+&\gamma \lbrack g(Y,Z)X-g(X,Z)Y],  \label{2.6}
\end{eqnarray}%
for all $X,Y,Z\in \chi (M)$, where $S(U,V)=g(QU,V)$ and $\alpha ,\beta
,\gamma $ are real constants. In particular, the quasi-conformal like
curvature tensor $\tilde{W}$ is

(i) the Riemann curvature tensor $R$ if $\alpha =\beta =\gamma =0,$

(ii) the conharmonic curvature tensor $H$ \cite{Ishii} if $\alpha =\beta =-%
\frac{1}{2{\Large n}-1},$ $\gamma =0,$

(iii) the conformal curvature tensor $C$ \cite{Eisenhart} if $\alpha =\beta
=-\frac{1}{2{\Large n}-1},$ $\gamma =\frac{r}{2n(2n-1)},$

(iv) the concircular curvature tensor $\tilde{C}$ \cite{Yano} if $\alpha
=\beta =0$, $\gamma =-\frac{r}{2n(2n+1)},$

(v) the projective curvature tensor $P$ \cite{Yano} if $\alpha =-\frac{1}{2%
{\Large n}},$ $\beta =\gamma =0$,

(vi) the ${\Large m}$-projective curvature tensor $\mathcal{M}$ \cite{PGP}
if $\alpha =\beta =$ $-\frac{1}{4{\Large n}},$ $\gamma =0,$

(vii) the $\mathcal{W}_{1}$-curvature tensor \cite{PGP1} if $\alpha =\frac{1%
}{n-1},$ $\beta =$ $0,$ $\gamma =0,$

(viii) the $\mathcal{W}_{2}$-curvature tensor \cite{PGP} if $\alpha =0,$ $%
\beta =$ $-\frac{1}{n-1},$ $\gamma =0,$

(ix) the $\mathcal{W}_{4}$-curvature tensor \cite{PGP2} if $\alpha =0,$ $%
\beta =0$, $\gamma =-\frac{1}{n-1}$.

The present paper deals with a study of 3-dimensional generalized ($\kappa
,\mu $)-contact metric manifolds with some semi-symmetric type of
geometric structure which cover 8 equivalent classes given by Shaikh and Kundu
\cite{Shaikh2} (see also \cite{kbmrb}, \cite{kbmb2}). Our paper is organized
as follows: Section 2 is concerned with some known results of
3-dimensional generalized ($\kappa ,\mu $)-contact metric manifolds.
Thereafter, in Section 3, we study 3-dimensional generalized ($%
\kappa ,\mu $)-contact metric manifolds with $\tilde{W}\cdot R=0$ and $%
\tilde{W}\cdot H=0.$ Here, we observed that a 3-dimensional non-Sasakian
generalized ($\kappa ,\mu $)-contact metric manifold satisfying the relation
$\tilde{W}\cdot R=0$ is either a ($\kappa ,\mu $)-contact metric manifold,
flat, or of constant $\xi $-sectional curvature $\kappa <1$ and constant $%
\phi $-sectional curvature $-\kappa $. It is also found that a $3$%
-dimensional non-Sasakian generalized ($\kappa ,\mu $)-contact metric
manifold belonging to the classes $G_{i}$ ($i=1,2,...,8$) is either a ($%
\kappa ,\mu $)-contact metric manifold, flat, or of constant $\xi $%
-sectional curvature $\kappa <1$ and constant $\phi $-sectional curvature $%
-\kappa $. Finally, in Section 4, we study a 3-dimensional generalized ($%
\kappa ,\mu $)-contact metric manifold with $\tilde{W}\cdot S=0$.

\section{Preliminaries}

In this section, we present some basic facts about contact metric manifolds.
We refer the reader to \cite{Blair} for more details. A differentiable
manifold $M^{2n+1}$ is called a contact manifold if it carries a global $1$%
-form $\eta $ such that $\eta \wedge (d\eta )^{2n+1}\neq 0$ everywhere on $%
M^{2n+1}.$ The form $\eta $ is usually called the contact form of $M^{2n+1}.$
It is well known that a contact metric manifold admits an almost contact
metric structure ($\phi ,\xi ,\eta ,g$), i.e. a global vector field $\xi $,
which is called the characteristic vector field, a $(1,1)$-tensor field $%
\phi $, and a Riemannian metric $g$ such that%
\begin{equation*}
\phi ^{2}X=-X+\eta (X)\xi,
\end{equation*}%
\begin{equation*}
\eta (\xi )=1,
\end{equation*}%
\begin{equation*}
g(\phi X,\phi Y)=g(X,Y)-\eta (X)\eta (Y),
\end{equation*}%
for all $X,Y$ $\in $ $\chi (M).$ Moreover, ($\phi ,\xi ,\eta ,g$) can be
chosen such that%
\begin{equation*}
d\eta (X,Y)=g(X,\phi Y),
\end{equation*}%
for all $X,Y$ $\in $ $\chi (M),$ and we call the structure a contact metric
structure. A manifold $M^{2n+1}$ carrying such a structure is said to be a
contact metric manifold and it is denoted by ($M^{2n+1},\phi ,\xi ,\eta ,g$%
). As a consequence of the above relations, we have%
\begin{equation*}
\eta (\xi )=1,
\end{equation*}%
\begin{equation*}
\phi \xi =0,
\end{equation*}%
\begin{equation*}
\eta \circ \phi =0 \ \text{and} \ d\eta (X,\xi )=0.
\end{equation*}%
If $\nabla $ denotes the Riemannian connection of ($M^{2n+1},\phi ,\xi ,\eta
,g$), then following \cite{blair}, we define the $(1,1)$-tensor fields $h$
and $l$ by%
\begin{equation*}
h=\frac{1}{2}(\pounds _{\xi }\phi )\text{ and }l=R(\cdot ,\xi )\xi ,
\end{equation*}%
where $\pounds _{\xi }$ is the Lie differentiation in the direction of $\xi $
and $R$ is the Riemann curvature tensor. The tensor fields $h$ and $l$ are
self-adjoint and they satisfy%
\begin{equation*}
h\xi =0,\text{ }l\xi =0,
\end{equation*}%
\begin{equation*}
trh=tr\phi h=0,\text{ }\phi h=-h\phi .
\end{equation*}%
Since $h$ anticommutes with $\phi $, if $X\neq 0$ is an eigenvector of $h$
corresponding to the eigenvalue $\theta $, then $\phi X$ is also an
eigenvector of $h$ corresponding to the eigenvalue $-\theta $. Therefore, on
any contact metric manifold ($M^{2n+1},\phi ,\xi ,\eta ,g$), the following
formulas are valid:%
\begin{equation*}
\nabla _{X}\xi =-\phi X-\phi hX\text{ (and so }\nabla _{\xi }\xi =0\text{),}
\end{equation*}%
\begin{equation*}
\nabla _{\xi }h=\phi -\phi l-\phi h^{2},
\end{equation*}%
\begin{equation*}
\phi l\phi -l=2(\phi ^{2}+h^{2}).
\end{equation*}%
A contact metric structure ($\phi ,\xi ,\eta ,g$) on $M^{2n+1}$ gives rise
to an almost complex structure on $M\times
%TCIMACRO{\U{211d} }%
%BeginExpansion
\mathbb{R}
%EndExpansion
.$ If this structure is integrable, then the contact metric manifold ($%
M^{2n+1},\phi ,\xi ,\eta ,g$) is said to be Sasakian. Equivalently, a
contact metric manifold ($M^{2n+1},\phi ,\xi ,\eta ,g$) is Sasakian if and
only if%
\begin{equation*}
R(X,Y)\xi =\eta (Y)X-\eta (X)Y,
\end{equation*}%
for all $X,Y$ $\in $ $\chi (M).$

By a generalized ($\kappa ,\mu $)-contact metric manifold, we mean a $3$-dimensional
contact metric manifold such that%
\begin{equation*}
R(X,Y)\xi =(\kappa I+\mu h)[\eta (Y)X-\eta (X)Y],
\end{equation*}%
for all $X,Y$ $\in $ $\chi (M)$, where $\kappa ,\mu $ are smooth nonconstant
real functions on $M^{2n+1}.$ In the special case where $\kappa ,\mu $ are
constant, then ($M^{2n+1},\phi ,\xi ,\eta ,g$) is called a ($\kappa ,\mu $%
)-contact metric manifold. We note that in any Sasakian manifold, $h=0$ and $\kappa =1.$

Furthermore, in a generalized ($\kappa ,\mu $)-contact metric manifold ($%
M^{2n+1},\phi ,\xi ,\eta ,g$), we have the following \cite{blair}, \cite%
{Koufogiorgos1}, \cite{Koufogiorgos2}:
\begin{eqnarray}
(\nabla _{X}h)Y &=&\{(1-\kappa )g(X,\phi Y)-g(X,\phi hY)\}\xi  \notag \\
&&-\eta (Y)\{(1-\kappa )\phi X+\phi hX\}-\mu \eta (X)\phi hY,  \label{2.1}
\end{eqnarray}%
\begin{equation}
(\nabla _{X}\phi )Y=\{g(X,Y)+g(X,hY)\}\xi -\eta (Y)(X+hX),  \label{2.2}
\end{equation}%
\begin{equation}
h^{2}=(\kappa -1)\phi ^{2}\text{, \ }\kappa =\frac{trl}{2}\leq 1,
\label{2.2A}
\end{equation}%
\begin{equation}
\xi \kappa =0,\text{ \ }\xi r=0,  \label{2.2B}
\end{equation}%
\begin{equation}
h\func{grad}\mu =\func{grad}\kappa .  \label{2.2C}
\end{equation}%
For a $3$-dimensional Riemannian manifold, the conformal curvature tensor $C$
vanishes. Thus, the Riemann curvature tensor $R$ and the Ricci tensor $S$
for a $3$-dimensional generalized ($\kappa ,\mu $)-contact metric manifold
are given by \cite{Shaikh}%
\begin{eqnarray}
R(U,V)Z &=&-(\kappa +\mu )\{g(V,Z)U-g(U,Z)V\}  \notag \\
&&+(2\kappa +\mu )\{g(V,Z)\eta (U)\xi -g(U,Z)\eta (V)\xi  \notag \\
&&+\eta (V)\eta (Z)U-\eta (U)\eta (Z)V\}+\mu \{g(V,Z)hU  \notag \\
&&-g(U,Z)hV+g(hV,Z)U-g(hU,Z)V\},  \label{2.3}
\end{eqnarray}%
\begin{equation}
S(U,V)=-\mu g(U,V)+\mu g(hU,V)+(2\kappa +\mu )\eta (U)\eta (V),  \label{2.4}
\end{equation}%
for all $U,V,Z\in \chi (M)$. We see that on a 3-dimensional generalized ($%
\kappa ,\mu $)-contact metric manifold, the scalar curvature $r$ equals to%
\begin{equation}
r=2(\kappa -\mu ).  \label{2.5}
\end{equation}
Now, using (\ref{2.4}) and (\ref{2.5}), (\ref{2.6}) yields%
\begin{eqnarray}
\tilde{W}(U,V)Z &=&\{-(\alpha +\beta +1)\mu -\kappa +\gamma
)\}[g(V,Z)U-g(U,Z)V]  \notag \\
&&+(\alpha +1)\mu \{g(hV,Z)U-g(hU,Z)V\}  \notag \\
&&+(\beta +1)\mu \{g(V,Z)hU-g(U,Z)hV\}  \notag \\
&&+(\alpha +1)(2\kappa +\mu )\{\eta (V)\eta (Z)U-\eta (U)\eta (Z)V\}  \notag
\\
&&+(\beta +1)(2\kappa +\mu )\{g(V,Z)\eta (U)-g(U,Z)\eta (V)\}\xi .
\label{2.7}
\end{eqnarray}%
From (\ref{2.7}), we have%
\begin{eqnarray}
\tilde{W}(\xi ,V)Z &=&\{-(\alpha +\beta +1)\mu -\kappa +\gamma \}[g(V,Z)\xi
-\eta (Z)V]  \notag \\
&&+(\alpha +1)\mu g(hV,Z)\xi -(\beta +1)\mu \eta (Z)hV  \notag \\
&&+(\alpha +1)(2\kappa +\mu )\{\eta (V)\eta (Z)\xi -\eta (Z)V\}  \notag \\
&&+(\beta +1)(2\kappa +\mu )\{g(V,Z)-\eta (Z)\eta (V)\}\xi ,  \label{2.8}
\end{eqnarray}%
\begin{eqnarray}
\tilde{W}(\xi ,V)\xi &=&\{(2\alpha +1)\kappa -\beta \mu +\gamma \}\{\eta
(V)\xi -V\}  \notag \\
&&-(\beta +1)\mu hV,  \label{2.9}
\end{eqnarray}%
\begin{equation}
S(U,\xi )=2\kappa \eta (U),  \label{2.10}
\end{equation}%
\begin{equation}
H(U,V)Z=(\mu -\kappa )\{g(V,Z)U-g(U,Z)V\},  \label{2.11}
\end{equation}%
\begin{equation}
H(\xi ,V)Z=(\mu -\kappa )\{g(V,Z)\xi -\eta (Z)V\}.  \label{2.12}
\end{equation}

\begin{lemma}
\label{l1} \cite{Koufogiorgos1} Let $M^{2n+1}$ be a non-Sasakian, generalized ($\kappa
,\mu $)-contact metric manifold. If $\kappa ,$ $\mu $\ satisfy the condition
$\pi \kappa +\rho \mu =\sigma $ ($\pi ,\rho ,\sigma $ are constants), then $%
\kappa ,$ $\mu $ are constant.
\end{lemma}

\begin{lemma}
\label{l2} \cite{Blair2} A $3$-dimensional contact metric manifold ($%
M^{3},\phi ,\xi ,\eta ,g$) with $\phi Q=Q\phi $ is either Sasakian, flat, or
of constant $\xi $-sectional curvature $\kappa <1$ and constant $\phi $%
-sectional curvature $-\kappa $.
\end{lemma}

\section{3-dimensional generalized ($\protect\kappa ,\protect\mu $)-contact
metric manifolds admitting some semi-symmetric structure}

In this section, we study $\tilde{W}\cdot R=0$ and $\tilde{W}\cdot H=0$ to
cover all the eight equivalent classes given by\cite{Shaikh2}. In \cite%
{Shaikh2} (also see \cite{kbmrb}, \cite{kbmb2}), the authors investigated
the equivalency of the various geometric structures. They have established
the following conditions:

i) $R\cdot H=0$ and $R\cdot C=0$ are equivalent and named such a class by $%
G_{1}$;

ii) $C\cdot C=0$ and $C\cdot H=0$ are equivalent and named such a class by $%
G_{2}$;

iii) $\tilde{C}\cdot C=0$ and $\tilde{C}\cdot H=0$ are equivalent and named
such a class by $G_{3}$;

iv) $H\cdot C=0$ and $H\cdot H=0$ are equivalent and named such a class by $%
G_{4}$;

v) $R\cdot R=0,$ $R\cdot P=0$, $R\cdot \tilde{C}=0,$ $R\cdot P^{\ast }=0,$
$R\cdot \mathcal{M}=0,$ $R\cdot \mathcal{W}_{i}=0$ and $R\cdot \mathcal{W}%
_{i}^{\ast }=0$ (for all $i=1,2,...,9$) are equivalent and named such a
class by $G_{5}$;

vi) $C\cdot R=0,$ $C\cdot P=0$, $C\cdot \tilde{C}=0,$ $C\cdot P^{\ast }=0,$ $%
C\cdot \mathcal{M}=0,$ $C\cdot \mathcal{W}_{i}=0$ and $C\cdot \mathcal{W}%
_{i}^{\ast }=0$ (for all $i=1,2,...,9$) are equivalent and named such a
class by $G_{6}$;

vii) $\tilde{C}\cdot R=0,$ $\tilde{C}\cdot P=0$, $\tilde{C}\cdot \tilde{C}=0$,
$\tilde{C}\cdot P^{\ast }=0,$ $\tilde{C}\cdot \mathcal{M}=0,$ $\tilde{C}%
\cdot \mathcal{W}_{i}=0$ and $\tilde{C}\cdot \mathcal{W}_{i}^{\ast }=0$ (for
all $i=1,2,...,9$) are equivalent and named such a class by $G_{7}$;

viii) $H\cdot R=0,$ $H\cdot P=0$, $H\cdot \tilde{C}=0,$ $H\cdot P^{\ast }=0,$
$H\cdot \mathcal{M}=0,$ $H\cdot \mathcal{W}_{i}=0$ and $H\cdot \mathcal{W}%
_{i}^{\ast }=0$ (for all $i=1,2,...,9$) are equivalent and named such a
class by $G_{8}$.

\subsection{3-dimensional generalized ($\protect\kappa ,\protect\mu $%
)-contact metric manifolds with $\tilde{W}\cdot R=0$}

Let ($M^{3},\phi ,\xi ,\eta ,g$) be a $3$-dimensional generalized ($\kappa
,\mu $)-contact metric manifold satisfying the relation $\tilde{W}\cdot R=0.$
Now,
\begin{eqnarray}
&&\tilde{W}(X,Y)R(U,V)Z-R(\tilde{W}(X,Y)U,V)Z  \notag \\
&&-R(U,\tilde{W}(X,Y)V)Z-R(U,V)\tilde{W}(X,Y)Z=0.  \label{3.1}
\end{eqnarray}%
Setting $X=U=\xi $ in (\ref{3.1}), we have
\begin{eqnarray}
&&\tilde{W}(\xi ,Y)R(\xi ,V)Z-R(\tilde{W}(\xi ,Y)\xi ,V)Z  \notag \\
&&-R(\xi ,\tilde{W}(\xi ,Y)V)Z-R(\xi ,V)\tilde{W}(\xi ,Y)Z=0.  \label{3.2}
\end{eqnarray}%
Now, by virtue of (\ref{2.3}) and (\ref{2.7}), (\ref{3.2}) yields
\begin{eqnarray}
&&\kappa \{g(V,Z)\tilde{W}(\xi ,Y)\xi -\eta (Z)\tilde{W}(\xi ,Y)V\}  \notag
\\
&&+\mu \{g(hV,Z)\tilde{W}(\xi ,Y)\xi -\eta (Z)\tilde{W}(\xi ,Y)hV\}-(\alpha
+1)\mu g(hY,Z)R(\xi ,V)\xi   \notag \\
&&-(2\kappa \alpha +\kappa +\gamma -\beta \mu )\{\eta (Y)R(\xi ,V)Z-R(Y,V)Z\}
\notag \\
&&+(\beta +1)\mu R(hY,V)Z+(\beta +1)\mu \eta (Z)R(\xi ,V)hY  \notag \\
&&+(2\kappa \alpha +\kappa +\gamma -\beta \mu )R(\xi ,Y)Z+(\beta +1)\mu \eta
(V)R(\xi ,hY)Z  \notag \\
&&+(\alpha \mu +\beta \mu +\kappa +\mu -\gamma )\{g(Y,Z)R(\xi ,V)\xi -\eta
(Z)R(\xi ,V)Y\}  \notag \\
&&-(\alpha +1)(2\kappa +\mu )\eta (Z)\{\eta (Y)R(\xi ,V)\xi -R(\xi ,V)Y\}
\notag \\
&&-(\beta +1)(2\kappa +\mu )\{g(Y,Z)-\eta (Y)\eta (Z)\}R(\xi ,V)\xi =0.
\label{3.3}
\end{eqnarray}%
Now, using (\ref{2.3}) and (\ref{2.8}) in (\ref{3.3}), we get
\begin{eqnarray}
&&\{\kappa g(V,Z)+\mu g(hV,Z)\}[\{(2\alpha +1)\kappa -\beta \mu +\gamma
\}\{\eta (Y)\xi -Y\}  \notag \\
&&-(\beta +1)\mu hY]-\kappa \eta (Z)[\{-(\alpha +\beta +1)\mu -\kappa
+\gamma )\}\{g(Y,V)\xi -\eta (V)Y\}  \notag \\
&&+(\alpha +1)\mu g(hY,V)\xi -(\beta +1)\mu \eta (V)hY  \notag \\
&&+(\alpha +1)(2\kappa +\mu )\{\eta (Y)\eta (V)\xi -\eta (V)Y\}  \notag \\
&&+(\beta +1)(2\kappa +\mu )\{g(Y,V)-\eta (V)\eta (Y)\}\xi ]  \notag \\
&&-\mu \eta (Z)[\{-(\alpha +\beta +1)\mu -\kappa +\gamma )\}g(Y,hV)\xi
\notag \\
&&+(\alpha +1)\mu g(hY,hV)\xi +(\beta +1)(2\kappa +\mu )g(Y,hV)]  \notag \\
&&+[(\alpha \mu +\beta \mu +\kappa +\mu -\gamma )g(Y,Z)-(\alpha +1)\mu
g(hY,Z)  \notag \\
&&-(\beta +1)(2\kappa +\mu )\{g(Y,Z)-\eta (Y)\eta (Z)\}  \notag \\
&&-(\alpha +1)(2\kappa +\mu )\eta (Z)\eta (Y)][\kappa \{\eta (V)\xi -V\}-\mu
hV]  \notag \\
&&+(2\kappa \alpha +\kappa +\gamma -\beta \mu )R(Y,V)Z+(\beta +1)\mu R(hY,V)Z
\notag \\
&&-(2\kappa \alpha +\kappa +\gamma -\beta \mu )\eta (Y)[\kappa \{g(V,Z)\xi
-\eta (Z)V\}  \notag \\
&&+\mu \{g(hV,Z)\xi -\eta (Z)hV\}]  \notag \\
&&+(\beta +1)\mu \eta (Z)[\kappa g(V,hY)\xi +\mu g(hV,hY)\xi ]  \notag \\
&&+(2\kappa \alpha +\kappa +\gamma -\beta \mu )[\kappa \{g(Y,Z)\xi -\eta
(Z)Y\}+\mu \{g(hY,Z)\xi -\eta (Z)hY\}]  \notag \\
&&+(\beta +1)\mu \eta (V)[\kappa \{g(hY,Z)\xi -\eta (Z)hY\}+\mu
\{g(h^{2}Y,Z)\xi -\eta (Z)h^{2}Y\}]  \notag \\
&&+\eta (Z)[(\alpha +1)(2\kappa +\mu )-(\alpha \mu +\beta \mu +\kappa +\mu
-\gamma )]  \notag \\
&&[\kappa \{g(V,Y)\xi -\eta (Y)V\}+\mu \{g(hV,Y)\xi -\eta (Y)hV\}]=0.
\label{3.4}
\end{eqnarray}%
Taking the inner product with $T,$ then contracting over $V$ and $Z$, we obtain
\begin{eqnarray}
&&\kappa \lbrack 3\{(2\alpha +1)\kappa -\beta \mu +\gamma \}+2(\alpha \mu
+\beta \mu +\kappa +\mu -\gamma )  \notag \\
&&-(\alpha +1)(2\kappa +\mu )+(2\kappa \alpha +\kappa +\gamma -\beta \mu )
\notag \\
&&-(\beta +1)(2\kappa +\mu )][\eta (Y)\eta (T)-g(Y,T)]  \notag \\
&&-3\kappa (\beta +1)\mu g(hY,T)-2\kappa (2\kappa \alpha +\kappa +\gamma
-\beta \mu )\eta (Y)\eta (T)  \notag \\
&&+(2\kappa \alpha +\kappa +\gamma -\beta \mu )S(Y,T)-(\beta +1)\mu S(hY,T)
\notag \\
&&-(2\kappa \alpha +\kappa +\gamma -\beta \mu )\mu g(hY,T)-(\beta +1)\mu
^{2}g(h^{2}Y,T)  \notag \\
&&-(\alpha \mu +\beta \mu +\kappa +\mu -\gamma )\mu g(Y,T)=0.  \label{3.5}
\end{eqnarray}%
That is, \
\begin{eqnarray}
&&\kappa \lbrack 3\{(2\alpha +1)\kappa -\beta \mu +\gamma \}+2(\alpha \mu
+\beta \mu +\kappa +\mu -\gamma )  \notag \\
&&-(\alpha +1)(2\kappa +\mu )+(2\kappa \alpha +\kappa +\gamma -\beta \mu )
\notag \\
&&-(\beta +1)(2\kappa +\mu )][\eta (Y)\xi -Y]  \notag \\
&&-3\kappa (\beta +1)\mu hY-2\kappa (2\kappa \alpha +\kappa +\gamma -\beta
\mu )\eta (Y)\xi   \notag \\
&&+(2\kappa \alpha +\kappa +\gamma -\beta \mu )QY-(\beta +1)\mu QhY  \notag
\\
&&-(2\kappa \alpha +\kappa +\gamma -\beta \mu )\mu hY-(\beta +1)\mu
^{2}h^{2}Y  \notag \\
&&-(\alpha \mu +\beta \mu +\kappa +\mu -\gamma )\mu Y=0.  \label{3.6}
\end{eqnarray}%
Operating $h$ and then taking the trace of (\ref{3.6}), we have%
\begin{equation}
(\beta +1)\mu (\mu -3\kappa )=0.  \label{3.7}
\end{equation}%
Therefore from (\ref{3.7}) we obtain either $\mu =0$ or $\mu -3\kappa =0$.
Now, for $\mu =0$ we conclude $Q\phi =\phi Q$ as we know that $Q\phi -\phi
Q=2\mu h\phi $ and for $\mu -3\kappa =0$ we get from Lemma \ref{l1} that $%
\kappa ,$ $\mu $ are constant. Hence from Lemma \ref{l1} and Lemma \ref{l2}
we can conclude the following:

\begin{theorem}
A $3$-dimensional non-Sasakian generalized ($\kappa ,\mu $)-contact metric
manifold satisfying the relation $\tilde{W}\cdot R=0$ is either a ($\kappa
,\mu $)-contact metric manifold, flat, or of constant $\xi $-sectional
curvature $\kappa <1$ and constant $\phi $-sectional curvature $-\kappa $.
\end{theorem}

\subsection{3-dimensional generalized ($\protect\kappa ,\protect\mu $%
)-contact metric manifolds with $\tilde{W}\cdot H=0$}

Let ($M^{3},\phi ,\xi ,\eta ,g$) be a $3$-dimensional generalized ($\kappa
,\mu $)-contact metric manifold satisfying the relation $\tilde{W}\cdot H=0.$
Now,
\begin{eqnarray}
&&\tilde{W}(X,Y)H(U,V)Z-H(\tilde{W}(X,Y)U,V)Z  \notag \\
&&-H(U,\tilde{W}(X,Y)V)Z-H(U,V)\tilde{W}(X,Y)Z=0.  \label{3.8}
\end{eqnarray}%
Setting $X=U=\xi $ in (\ref{3.8}), we have
\begin{eqnarray}
&&\tilde{W}(\xi ,Y)H(\xi ,V)Z-H(\tilde{W}(\xi ,Y)\xi ,V)Z  \notag \\
&&-H(\xi ,\tilde{W}(\xi ,Y)V)Z-H(\xi ,V)\tilde{W}(\xi ,Y)Z=0.  \label{3.9}
\end{eqnarray}%
Now, by virtue of (\ref{2.11}) and (\ref{2.7}), (\ref{3.9}) yields
\begin{eqnarray}
&&\kappa \{g(V,Z)\tilde{W}(\xi ,Y)\xi -\eta (Z)\tilde{W}(\xi ,Y)V\}  \notag
\\
&&+\mu \{g(hV,Z)\tilde{W}(\xi ,Y)\xi -\eta (Z)\tilde{W}(\xi ,Y)hV\}-(\alpha
+1)\mu g(hY,Z)H(\xi ,V)\xi  \notag \\
&&-(2\kappa \alpha +\kappa +\gamma -\beta \mu )\{\eta (Y)H(\xi ,V)Z-H(Y,V)Z\}
\notag \\
&&+(\beta +1)\mu H(hY,V)Z+(\beta +1)\mu \eta (Z)H(\xi ,V)hY  \notag \\
&&+(2\kappa \alpha +\kappa +\gamma -\beta \mu )H(\xi ,Y)Z+(\beta +1)\mu \eta
(V)H(\xi ,hY)Z  \notag \\
&&+(\alpha \mu +\beta \mu +\kappa +\mu -\gamma )\{g(Y,Z)H(\xi ,V)\xi -\eta
(Z)H(\xi ,V)Y\}  \notag \\
&&-(\alpha +1)(2\kappa +\mu )\eta (Z)\{\eta (Y)H(\xi ,V)\xi -H(\xi ,V)Y\}
\notag \\
&&-(\beta +1)(2\kappa +\mu )\{g(Y,Z)-\eta (Y)\eta (Z)\}H(\xi ,V)\xi =0.
\label{3.10}
\end{eqnarray}%
Now, using (\ref{2.7}) and (\ref{2.11}) in (\ref{3.10}), we get
\begin{eqnarray}
&&\{\kappa g(V,Z)+\mu g(hV,Z)\}[\{(2\alpha +1)\kappa -\beta \mu +\gamma
\}\{\eta (Y)\xi -Y\}  \notag \\
&&-(\beta +1)\mu hY]-\kappa \eta (Z)[\{-(\alpha +\beta +1)\mu -\kappa
+\gamma )\}\{g(Y,V)\xi -\eta (V)Y\}  \notag \\
&&+(\alpha +1)\mu g(hY,V)\xi -(\beta +1)\mu \eta (V)hY  \notag \\
&&+(\alpha +1)(2\kappa +\mu )\{\eta (Y)\eta (V)\xi -\eta (V)Y\}  \notag \\
&&+(\beta +1)(2\kappa +\mu )\{g(Y,V)-\eta (V)\eta (Y)\}\xi ]  \notag \\
&&-\mu \eta (Z)[\{-(\alpha +\beta +1)\mu -\kappa +\gamma )\}g(Y,hV)\xi
\notag \\
&&+(\alpha +1)\mu g(hY,hV)\xi +(\beta +1)(2\kappa +\mu )g(Y,hV)]  \notag \\
&&+[(\alpha \mu +\beta \mu +\kappa +\mu -\gamma )g(Y,Z)-(\alpha +1)\mu
g(hY,Z)  \notag \\
&&-(\beta +1)(2\kappa +\mu )\{g(Y,Z)-\eta (Y)\eta (Z)\}  \notag \\
&&-(\alpha +1)(2\kappa +\mu )\eta (Z)\eta (Y)][(\mu -\kappa )\{\eta (V)\xi
-V\}]  \notag \\
&&+(2\kappa \alpha +\kappa +\gamma -\beta \mu )H(Y,V)Z+(\beta +1)\mu H(hY,V)Z
\notag \\
&&-(2\kappa \alpha +\kappa +\gamma -\beta \mu )\eta (Y)(\mu -\kappa
)[\{g(V,Z)\xi -\eta (Z)V\}  \notag \\
&&+\{g(Y,Z)\xi -\eta (Z)Y\}]+\mu (\beta +1)(\mu -\kappa )\eta (Z)g(V,hY)\xi
\notag \\
&&+\mu (\beta +1)(\mu -\kappa )\eta (V)\{g(hY,Z)\xi -\eta (Z)hY\}  \notag \\
&&-(\alpha \mu +\beta \mu +\kappa +\mu -\gamma )(\mu -\kappa )\eta
(Z)\{g(Y,V)\xi -\eta (Y)V\}  \notag \\
&&+\eta (Z)(\alpha +1)(2\kappa +\mu )(\mu -\kappa )[\{g(V,Y)\xi -\eta (Y)V\}=0.  \label{3.11}
\end{eqnarray}%
Contracting over $V$ and $Z$, we obtain%
\begin{eqnarray}
&&3\kappa \lbrack \{(2\alpha +1)\kappa -\beta \mu +\gamma \}\{\eta (Y)\xi
-Y\}  \notag \\
&&-(\beta +1)\mu hY]-\kappa \lbrack \{-(\alpha +\beta +1)\mu -\kappa +\gamma
)\}\{\eta (Y)\xi -Y\}  \notag \\
&&-(\beta +1)\mu hY+(\alpha +1)(2\kappa +\mu )\{\eta (Y)\xi -Y\}]  \notag \\
&&-(\mu -\kappa )[-(\alpha \mu +\beta \mu +\kappa +\mu -\gamma )\{\eta
(Y)\xi +Y\}+(\alpha +1)(2\kappa +\mu )\eta (Y)\xi  \notag \\
&&-(\alpha +1)\mu hY-(\alpha +1)(2\kappa +\mu )\eta (Y)\xi -(\beta
+1)(2\kappa +\mu )\{Y-\eta (Y)\xi \}]  \notag \\
&&-(2\kappa \alpha +\kappa +\gamma -\beta \mu )(\mu -\kappa )\{\eta (Y)\xi
-Y\}+(\beta +1)\mu (\mu -\kappa )hY =0.  \label{3.12}
\end{eqnarray}%
Operating $h$ and the taking the trace of (\ref{3.12}) we have%
\begin{equation}
\mu \lbrack (\alpha +\beta +2)\mu -(\alpha +3\beta +4)\kappa ]=0.
\label{3.13}
\end{equation}%
Therefore from (\ref{3.13}) we obtain either $\mu =0$ or $(\alpha +\beta
+2)\mu -(\alpha +3\beta +4)\kappa =0$. Now, for $\mu =0$ we conclude $Q\phi
=\phi Q$ as we know that $Q\phi -\phi Q=2\mu h\phi $ and for $(\alpha +\beta
+2)\mu -(\alpha +3\beta +4)\kappa =0$ we get from Lemma \ref{l1} that $%
\kappa ,$ $\mu $ are constant. Hence from Lemma \ref{l1} and Lemma \ref{l2}
we can conclude the following:

\begin{theorem}
A $3$-dimensional non-Sasakian generalized ($\kappa ,\mu $)-contact metric
manifold satisfying the relation $\tilde{W}\cdot H=0$ is either a ($\kappa
,\mu $)-contact metric manifold, flat, or of constant $\xi $-sectional
curvature $\kappa <1$ and constant $\phi $-sectional curvature $-\kappa $.
\end{theorem}

Again, a quasi-conformal like curvature tensor $\tilde{W}$ reduces to the
Riemann curvature tensor $R$ for $\alpha =\beta =\gamma =0,$ to the
conharmonic curvature tensor $H$ for $\alpha =\beta =-\frac{1}{2{\Large n}-1}%
,$ $\gamma =0,$ to the conformal curvature tensor $C$ for $\alpha =\beta =-%
\frac{1}{2{\Large n}-1},$ $\gamma =\frac{r}{2n(2n-1)},$ and to the
concircular curvature tensor $\tilde{C}$ for $\alpha =\beta =0$, $\gamma =-%
\frac{r}{2n(2n+1)}.$ Thus the relation $\tilde{W}\cdot R=0$ reduces to $%
R\cdot R=0$ for $\alpha =\beta =\gamma =0,$ to $H\cdot R=0$ for $\alpha
=\beta =-\frac{1}{2{\Large n}-1},$ $\gamma =0,$ to $C\cdot R=0$ for $\alpha
=\beta =-\frac{1}{2{\Large n}-1},$ $\gamma =\frac{r}{2n(2n-1)}$ and to $%
\tilde{C}\cdot R=0$ for $\alpha =\beta =0$, $\gamma =-\frac{r}{2n(2n+1)}$
and $\tilde{W}\cdot H=0$ reduces to $R\cdot H=0$ for $\alpha =\beta =\gamma
=0,$ to $H\cdot H=0$ for $\alpha =\beta =-\frac{1}{2{\Large n}-1},$ $\gamma
=0,$ to $C\cdot H=0$ for $\alpha =\beta =-\frac{1}{2{\Large n}-1},$ $\gamma =%
\frac{r}{2n(2n-1)}$ and to $\tilde{C}\cdot H=0$ for $\alpha =\beta =0$, $%
\gamma =-\frac{r}{2n(2n+1)}.$

From the above discussion, we can conclude

\begin{theorem}
A $3$-dimensional non-Sasakian generalized ($\kappa ,\mu $)-contact metric
manifold belonging to the classes $G_{i}$ ($i=1,2,...,8$) is either a ($%
\kappa ,\mu $)-contact metric manifold, flat, or of constant $\xi $%
-sectional curvature $\kappa <1$ and constant $\phi $-sectional curvature $%
-\kappa $.
\end{theorem}

\section{3-dimensional generalized ($\protect\kappa ,\protect\mu $)-contact
metric manifolds with $\tilde{W}\cdot S=0$}

Let ($M^{3},\phi ,\xi ,\eta ,g$) be a $3$-dimensional generalized ($\kappa
,\mu $)-contact metric manifold satisfying the relation $\tilde{W}\cdot S=0.$
Now,
\begin{equation}
S(\tilde{W}(\xi ,Y)U,V)+S(U,\tilde{W}(\xi ,Y)V)=0.  \label{4.1}
\end{equation}%
Now, using (\ref{2.8}) in (\ref{4.1}) yields%
\begin{eqnarray}
&&\{-(\alpha +\beta +1)\mu -\kappa +\gamma )\}[2\kappa \eta (V)g(Y,U)-\eta
(U)S(V,Y)  \notag \\
&&+2\kappa \eta (U)g(Y,V)-\eta (V)S(U,Y)]  \notag \\
&&+2(\alpha +1)\kappa \mu \eta (V)g(hY,U)-(\beta +1)\mu \eta (U)S(V,hY)
\notag \\
&&+2(\alpha +1)\kappa \mu \eta (U)g(hY,V)-(\beta +1)\mu \eta (V)S(U,hY)
\notag \\
&&+(\alpha +1)(2\kappa +\mu )\eta (U)[2\kappa \eta (V)\eta (Y)-S(V,Y)]
\notag \\
&&+(\alpha +1)(2\kappa +\mu )\eta (V)[2\kappa \eta (U)\eta (Y)-S(U,Y)]
\notag \\
&&+2\kappa (\beta +1)(2\kappa +\mu )\{g(Y,U)-\eta (U)\eta (Y)\}\eta (V)
\notag \\
&&+2\kappa (\beta +1)(2\kappa +\mu )\{g(Y,V)-\eta (V)\eta (Y)\}\eta (U)=0.
\label{4.2}
\end{eqnarray}%
Setting $U=\xi $ in (\ref{4.2}) yields%
\begin{eqnarray}
&&\{(\alpha +\beta +1)\mu +\kappa -\gamma )\}[-\mu g(Y,V)+\mu g(hY,V)  \notag
\\
&&+(2\kappa +\mu )\eta (Y)\eta (V)-2\kappa g(Y,V)]+2(\alpha +1)\kappa \mu
g(hY,V)  \notag \\
&&-(\beta +1)\mu ^{2}[g(h^{2}Y,V)-g(hY,V)]  \notag \\
&&+(\alpha +1)(2\kappa +\mu )[2\kappa \eta (V)\eta (Y)  \notag \\
&&+\mu g(Y,V)-\mu g(hY,V)-(2\kappa +\mu )\eta (Y)\eta (V)]  \notag \\
&&+2\kappa (\beta +1)(2\kappa +\mu )\{g(Y,V)-\eta (V)\eta (Y)\}=0.
\label{4.3}
\end{eqnarray}%
Notice that (\ref{4.3}) can be written as%
\begin{eqnarray}
&&\{(\alpha +\beta +1)\mu +\kappa -\gamma )\}[\mu hY-\mu Y  \notag \\
&&+(2\kappa +\mu )\eta (Y)\xi -2\kappa gY]+2(\alpha +1)\kappa \mu hY  \notag
\\
&&-(\beta +1)\mu ^{2}[(\kappa -1)\{-Y+\eta (Y)\xi \}-hY]  \notag \\
&&+(\alpha +1)(2\kappa +\mu )[2\kappa \eta (Y)\xi  \notag \\
&&+\mu Y-\mu hY-(2\kappa +\mu )\eta (Y)\xi ]  \notag \\
&&+2\kappa (\beta +1)(2\kappa +\mu )\{Y-\eta (Y)\xi \}=0.  \label{4.4}
\end{eqnarray}%
Operating $h$ and then taking the trace of (\ref{4.4}), we obtain%
\begin{equation}
[(2\beta +1)\mu +\kappa -\gamma ]\mu =0.  \label{4.5}
\end{equation}%
Therefore from (\ref{4.5}) we obtain either $\mu =0$ or $(2\beta +1)\mu
+\kappa =\gamma $. Now, for $\mu =0$ we conclude $Q\phi =\phi Q$ as we know
that $Q\phi -\phi Q=2\mu h\phi $ and for $(2\beta +1)\mu +\kappa =\gamma $
we get from Lemma \ref{l1} that $\kappa ,$ $\mu $ are constant. Hence from
Lemma \ref{l1} and Lemma \ref{l2} we can conclude the following:

\begin{theorem}
A $3$-dimensional non-Sasakian generalized ($\kappa ,\mu $)-contact metric
manifold satisfying the relation $\tilde{W}\cdot S=0$ is either a ($\kappa
,\mu $)-contact metric manifold, flat, or of constant $\xi $-sectional
curvature $\kappa <1$ and constant $\phi $-sectional curvature $-\kappa $.
\end{theorem}

A quasi-conformal like curvature tensor $\tilde{W}$ reduces to the Riemann
curvature tensor $R$ for $\alpha =\beta =\gamma =0,$ to the conharmonic
curvature tensor $H$ for $\alpha =\beta =-\frac{1}{2{\Large n}-1},$ $\gamma
=0,$ to the conformal curvature tensor $C$ for $\alpha =\beta =-\frac{1}{2%
{\Large n}-1},$ $\gamma =\frac{r}{2n(2n-1)},$ to the concircular curvature
tensor $\tilde{C}$ for $\alpha =\beta =0$, $\gamma =-\frac{r}{2n(2n+1)},$ to
the projective curvature tensor $P$ if $\alpha =-\frac{1}{2{\Large n}},$ $%
\beta =\gamma =0$, to the ${\Large m}$-projective curvature tensor $\mathcal{%
M}$ if $\alpha =\beta =$ $-\frac{1}{4{\Large n}},$ $\gamma =0,$ to the $%
\mathcal{W}_{1}$-curvature tensor if $\alpha =\frac{1}{n-1},$ $\beta =$ $0,$
$\gamma =0,$ to the $\mathcal{W}_{2}$-curvature tensor if $\alpha =0,$ $%
\beta =$ $-\frac{1}{n-1},$ $\gamma =0,$ and to the $\mathcal{W}_{4}$%
-curvature tensor if $\alpha =0,$ $\beta =0$ $,$ $\gamma =-\frac{1}{n-1}$.
Thus the relation $\tilde{W}\cdot S=0$ reduces to $R\cdot S=0$ for $\alpha
=\beta =\gamma =0,$ to $H\cdot S=0$ for $\alpha =\beta =-\frac{1}{2{\Large n}%
-1},$ $\gamma =0,$ to $C\cdot S=0$ for $\alpha =\beta =-\frac{1}{2{\Large n}%
-1},$ $\gamma =\frac{r}{2n(2n-1)}$, to $\tilde{C}\cdot S=0$ for $\alpha
=\beta =0$, $\gamma =-\frac{r}{2n(2n+1)},$ to $P\cdot S=0$ for $\alpha =-%
\frac{1}{2{\Large n}},$ $\beta =\gamma =0$, to $\mathcal{M}\cdot S=0$ for $%
\alpha =\beta =$ $-\frac{1}{4{\Large n}},$ $\gamma =0,$ to $\mathcal{W}%
_{1}\cdot S=0$ for $\alpha =\frac{1}{n-1},$ $\beta =$ $0,$ $\gamma =0,$ to $%
\mathcal{W}_{2}\cdot S=0$ for $\alpha =0,$ $\beta =$ $-\frac{1}{n-1},$ $%
\gamma =0,$ and to $\mathcal{W}_{4}\cdot S=0$ for $\alpha =0,$ $\beta =0$,
$\gamma =-\frac{1}{n-1}$.

From the above discussion, we can conclude

\begin{theorem}
A $3$-dimensional non-Sasakian generalized ($\kappa ,\mu $)-contact metric
manifold satisfying each of the relations $R\cdot S=0,$ $H\cdot S=0,$ $%
\tilde{C}\cdot S=0,$ $P\cdot S=0,$ $\mathcal{M}\cdot S=0,$ $\mathcal{W}%
_{1}\cdot S=0,$ $\mathcal{W}_{2}\cdot S=0$ or $\mathcal{W}_{4}\cdot S=0$
is either a ($\kappa ,\mu $)-contact metric manifold, flat, or of constant $%
\xi $-sectional curvature $\kappa <1$ and constant $\phi $-sectional
curvature $-\kappa $.
\end{theorem}

{\small {$^{\text{ }1}$Department of Mathematics, Raiganj University, Uttar
Dinajpur, India, } }

{\small E-mail address: raybakshimanoj@gmail.com.}

{\small {$^{\text{ }2}$Department of Mathematics, Kurseong College,
Kurseong, Darjeeling India, } }

{\small E-mail address: kanakkanti.kc@gmail.com. }

\end{document}